\documentstyle{amsart}

\title{Heitmann's Proof of the Direct Summand Conjecture in Dimension 3}

\author{Paul Roberts}

\newtheorem{thm}{Theorem}

\newtheorem{lem}{Lemma}

\newcommand{\m}{{\frak m}}

\begin{document}
\pagenumbering{roman}
\maketitle
\pagenumbering{arabic}

We describe the main ideas of Ray Heitmann's proof of the 
Direct Summand Conjecture in dimension 3
for a ring of mixed characteristic \cite{h}.  In the first section we describe the main
methods which are used and prove the necessary lemmas.  In the second section
we prove the main result of Heitmann's paper.  Finally, in the third section we
give a proof of the Canonical Element Conjecture using this result.

\section{An outline of the methods used in the proof}

Let $R$ be a Noetherian local ring of mixed characteristic $p$ of dimension 3.  We assume that
$R$ is an integrally closed domain  
and is 
a homomorphic image of a regular local ring.  In fact, the questions we
consider can be reduced to the case of a complete
integrally closed domain, and complete local
rings are homomorphic images of regular local rings
by the Cohen structure theorems. 

Let $(p,x,y)$ be a system of parameters for $R$.  Let $z$ be an
element of $R$ such that $p^Nz\in (x,y)$.  The aim is to find a finite extension of $R$ in
which $z$ is in the ideal generated by $x$ and $y$.  If this is
always possible, then it is not
difficult to prove the Direct Summand
Conjecture (or any of several equivalent conjectures).  This
aim is not quite reached, but a weaker result, which is also enough to prove these conjectures,
is proven instead.  However, to motivate the 
construction in the proof, we first discuss the question
of attempting to find an extension in which $z\in (x,y)$.  We assume throughout that $R$ is an
integrally closed domain.

The equation for which we are trying to find a solution is
$$z=wx+vy,$$
where $w$ and $v$ are in a finite extension of $R$.  This means that $w$ is integral over $R$,
so it is a root of a polynomial 
$$f(T)=T^n+a_1T^{n-1}+\cdots +a_n$$
with coefficients in $R$.  Since $z,x,$ and $y$ are in $R$, we can solve
the equation $z=wx+vy$ for $v$,
and $v$ will be in the field generated by $w$ over the quotient field of $R$:
$$v={{z-wx}\over {y}}.$$
From this equation we can use Taylor's formula 
to find a polynomial $h(T)$ whose root is $v$, as we show
below.  If we can find a polynomial 
$f(T)$ such that the corresponding polynomial $h(T)$ with
root
$v$ also has coefficients in $R$, we have solved the problem, since then both $w$ and $v$ are
integral over
$R$.

Suppose that $f(T)=T^n+a_1T^{n-1}+\cdots +a_n$ is a polynomial such that $f(w)=0$. We show how
to find the polynomial $h(T)$ with root $v$.  
Let $F(S,U)$ be $f(T)$ made homogeneous; that is, if
$f(T)=T^n+a_1T^{n-1}+\cdots +a_n$, we let
$$F(S,U)=S^n+a_1S^{n-1}U+\cdots + a_nU^n.$$
By Taylor's Theorem, we have
$$F(S-T,U)=\sum_{i=0}^n(-1)^i{1\over {i!}}F^{(i)}(S,U)T^i,$$
where $F^{(i)}(S,U)$ denotes the $i$th derivative of $F(S,T)$ with respect to the first
variable.  We will use this notation for the derivative of a homogeneous polynomial
with respect to the first variable, and we will refer to it as the derivative
of the polynomial, throughout the paper.  Applying this equation with
$S=z,T=z-wx,$ and $U=x$ we obtain
$$F(z-(z-wx),x)=\sum_{i=0}^n(-1)^i{1\over {i!}}F^{(i)}(z,x)(z-wx)^i.$$
On the other hand, we have
$$F(z-(z-wx),x)=F(wx,x)=x^nF(w,1)=x^nf(w)=0.$$
Hence if we let $$b_i=(-1)^{n-i}{1\over {(n-i)!}}F^{(n-i)}(z,x)$$ and let 
$$g(T)=T^n+b_1T^{n-1}+\cdots+b_n,$$
then $g(z-wx)=0$.  Finally, if we let $c_i=b_i/y^i$ and $h(T)=T^n+c_1T^{n-1}+\cdots +c_n$,
then $h((z-wx)/y)=0$.  

To summarize, we wish to find a polynomial $f(T)$ 
with coefficients in $R$ so that the resulting polynomial $h(T)$
also has coefficients in $R$.  
Tracing back through the construction of $h(T)$, we see
that this means that 
			$${1\over {(n-i)!}}F^{(n-i)}(z,x)\in y^iR$$
for $i=0,\ldots,n$.  

The strategy is to start with ${1\over {(n-1)!}}F^{(n-1)}(z,x)$, which is the next to the last
term in the Taylor expansion, and construct the previous terms inductively, essentially
by integration, finally ending up with a polynomial of degree $n$ with the required properties.

Before continuing the discussion, we recall some facts about the coefficients which arise in
the Taylor expansion and introduce some notation.  Let $f(T)=\sum_{k=0}^na_kT^{n-k}$ as above.
Then $${1\over {(n-i)!}}f^{(n-i)}(T)=\sum_{j=0}^i{{n-j}\choose{i-j}}a_jT^{i-j}.\eqno(1)$$
To see this, simply compute both sides for 
$f(T)=a_jT^{n-j}$: in computing the left hand
side we obtain
$${{(n-j)(n-j-1)\cdots (n-j-(n-i)+1)}\over{(n-i)!}}a_jT^{n-j-(n-i)}=$$
$$={{n-j}\choose {n-i}}a_jT^{i-j}
={{n-j}\choose {i-j}}a_jT^{i-j}.$$  
While this fact is quite straightforward, we mention
it explicitly to explain our interest in binomial coefficients of the
form $\displaystyle{{n-j}\choose {i-j}}$, for example in Lemma \ref{bincoef}.

Applying  equation (1) to the polynomial $F(S,U)$
evaluated at $z,x$, we have
  $${1\over {(n-i)!}}F^{(n-i)}(z,x)=\sum_{j=0}^i{{n-j}\choose{i-j}}a_jz^{i-j}x^j.$$

Since the main step in the construction involves integration
and adjusting constants, we introduce special notation for it.
Given a homogeneous polynomial $G(S,U)=\sum_{j=0}^i a_jS^{i-j}U^j$ of degree $i$ and an integer
$n>i$, we let
$$\mbox{Int}_n(G)(S,U)=(n-i)\sum_{j=0}^i a_j\left({{S^{i-j+1}\over{i-j+1}}}\right)U^j.$$
This definition depends not only on $G(S,U)$,
but also on $n$;
in fact, we are thinking of $G(S,U)$ as the $(n-i)$th term in the Taylor expansion
of a polynomial of degree $n$, and Int$_n(G)(S,U)$ is then a candidate for the previous term;
it is in fact the unique choice with coefficient of $U^{i+1}$ equal to zero.
That is the reason for including the factor $n-i$.

We now return to the main discussion.  We have a system of parameters $p,x,y$ for $R$
and an element $z$ with $p^Nz\in (x,y)$ for some integer $N$. Thus there exists an element $a$
of $R$ such that $p^Nz+ax\in yR$.   We want to find a homogeneous polynomial $F(S,U)$ of degree
$n=p^L$ for some
$L$ such that 
$${1\over {(n-i)!}}F^{(n-i)}(z,x)\in y^iR$$
for $i=0,\ldots,n$.  Since $n=p^L$, for $i=1$ this expression becomes
$${1\over {(n-1)!}}F^{(n-1)}(z,x)={{p^L}\choose 1}z+{{p^L-1}\choose 0}a_1x=p^Lz+a_1x.$$
From the above  condition on $z$, if we let $L=N$ and $a_1=a$, we have 
$$p^Lz+a_1x\in yR,$$
which is the required condition for $i=1$. We let $F_1(S,U)=p^LS+a_1U.$

As stated above, we wish to construct $F(S,U)$ by repeated integration.  If we knew that
Int$_{p^L}(F_1)(z,x)$ were in the ideal $(x^2,y^2)$, we could find an element $a_2$ such
that
$$\mbox{Int}_{p^L}(F_1)(z,x)+a_2x^2={{p^L}\choose 2}z^2+{{p^L-1}\choose 1}a_1zx+{{p^L-2}\choose
0}a_2x^2\in y^2R.$$
and continue the construction.  Unfortunately, this is not necessarily possible (we recall that
although we are discussing the question of trying to show that $z\in (x,y)R'$ for a finite
extension $R'$ of $R$, the actual result we prove is somewhat weaker).  However, the next lemma
shows that although we cannot solve this problem in $R$, we can solve it in $R[p^{-1}]$.

\begin{lem}\label{inxy} Let $\tilde R$ be a ring which contains 
a field of characteristic zero, and let $H(S,T)$
be a homogeneous polynomial of degree $n$ with coefficients in ${\tilde R}$.  Suppose 
that $x,y,z$ are
elements of ${\tilde R}$ such that $z\in(x,y)$ and such that
$$H^{(n-i)}(z,x)\in y^i{\tilde R}\;\;\mbox{for}\;\;i=0,\ldots,n-1.$$
Then $$H(z,x)\in (x^n,y^n).$$
\end{lem}
{\bf Proof.}
Since $z\in (x,y)$, we can write  $z=cx-dy$ for some $c$ and $d$ in ${\tilde R}$.  
Using Taylor's formula, we have
$$H(S+T,U)=\sum_{i=0}^n{1\over{(n-i)!}}H^{(n-i)}(S,U)T^{n-i}.$$
We apply this formula with $S=z$, $U=x$, and $T=dy$, giving
$$H(z+dy,x)=\sum_{i=0}^n{1\over{(n-i)!}}H^{(n-i)}(z,x)(dy)^{n-i}.$$ 
Since $H^{(n-i)}(z,x)\in y^i{\tilde R}$ for $i=0,\ldots,n-1$, each term in the sum on the right is
in $y^n{\tilde R}$ except possibly the term with $i=n$, which is $H(z,x)$.  Hence we can write
$$H(z+dy,x)=H(z,x)+ay^n$$
for some $a\in {\tilde R}$.  Since $z=cx-dy$, we also have
$$H(z+dy,x)=H(cx,x)=x^nH(c,1)\in x^n{\tilde R}.$$
Thus
$$H(z,x)=H(z+dy,x)-ay^n\in (x^n,y^n).$$

\bigskip

Let $u_2=\mbox{Int}_{p^L}(F_1)(z,x)$ as above.  Then  Lemma \ref{inxy}
 implies that $u_2\in(x^2,y^2)R[p^{-1}]$, 
so there is an integer $D$ such that $p^Du_2\in (x^2,y^2).$   
The main part of the proof consists of adjusting $u_2$ 
so that it lies in $(x^2,y^2)$ (and similarly adjusting the
corresponding element at the $i$th stage so that it lies in $(x^i,y^i)$.  There are
essentially three possibilities:
\begin{enumerate}
\item The element $u_2$ itself lies in $(x^2,y^2)$.  In this case, as outlined above, there is no problem.
\item There is an element $c\in R$ such that $u_2\equiv cxyz\;\; \mbox{modulo}
(x^2,y^2)$.  We recall that we had $p^Lz\in (x,y)$,
so $p^Lxyz\in (x^2,y^2)$.  In this case, we can add $(-c)xyz$ to $u_2$ and proceed as before (we will explain this in detail below).
\item The first two possibilities do not apply.  In this case we
have to multiply $u_2$ by a power
of $p$.  We recall that $L$ was chosen so that $p^Lz\in (x,y)$.  
Now we are raising the power
of $p$ needed, and we have to increase $L$.  Thus we are now looking for a polynomial of larger 
degree.  
\end{enumerate}

These three procedures form the basis of the construction.  One of the main facts which makes 
this construction
work is that after a certain point the third case will no longer occur.  The lemma we use here
is the following.

\begin{lem}\label{finitelh}
  Let $R$ be a ring as above.  For each $i$, let 
$$J_i=\{r\in R|p^Nr\in (x^i,y^i)\;\; 
\mbox{for\;\;some}\;N\}.$$
  Let $Q_i=J_i/(x^i,y^i)$.  Then
\begin{enumerate}
\item $Q_i$ has finite length for all $i$.
\item There exists an integer $n$ such that for all $k\ge n$, the map from $Q_n$ 
to $Q_{k}$ induced by multiplication by $x^{k-n}y^{k-n}$
is an isomorphism.
\end{enumerate}
\end{lem}
{\bf Proof.}
The fact that $Q_i$ has finite length follows immediately from the fact that $p,x,y$ is
a system of parameters.   

Since $x,y$ form part of a system of parameters and $R$ is assumed
to be integrally closed, $x,y$ form a regular sequence.  It follows
that the map from $Q_n$ to $Q_k$ induced by multiplication
by $x^{k-n}y^{k-n}$ is always injective.  To prove the second
statement it must be shown that this map is also surjective for
large enough $n$.

To prove this assertion
we use the fact that the local cohomology module $H^2_{\m}(R)$
has finite length, which follows from the assumptions we have made
on $R$.  More specifically, since $R$ is an integrally closed domain of
dimension three, its non-Cohen-Macaulay locus is supported at the
maximal ideal. It then follows from local duality and the assumption that
$R$ is a homomorphic image of a regular local ring that the local
cohomology module $H^2_{\m}(R)$ has finite length (see \cite{b-h},
section 3.5). We remark that this is the only place where
the assumption that $R$ is a homomorphic image 
of a regular local ring is used.  

The local cohomology module can be computed as the homology of the sequence
$$R_p\times R_x\times R_y\stackrel{d_1}{\to} R_{px}\times R_{py}\times R_{xy}
\stackrel{d_2}{\to} R_{pxy},$$
where the first map sends $(\alpha_1,\alpha_2,\alpha_3)$ to $(\alpha_1-\alpha_2,
\alpha_3-\alpha_1,\alpha_2-\alpha_3)$, and the second map sends 
$(\beta_1,\beta_2,\beta_3)$ to $\beta_1+\beta_2+\beta_3$.  Since $H^2_{\m}(R)$ is
finitely generated, there is an integer $n$ such that $H^2_{\m}(R)$ is generated
by elements which can be written in the form 
$$\left({{a}\over {(px)^n}},{{b}\over {(py)^n}},{{c}\over {(xy)^n}}\right).$$

Let $k\ge n$.  Let $r$ be an element of $R$ and $N$ an integer such that
$$p^Nr=sx^{k}+ty^{k}$$
for $s$ and $t$ in $R$.
We then have that the element 
$$\left({{-t}\over {p^Nx^k}},{{-s}\over {p^Ny^k}},{{r}\over {x^ky^k}}\right)\in R_{px}\times R_{py}
\times R_{xy}$$
is in the kernel of $d_2$.  By our choice of $n$, there is an element of the form
$$\left({{a}\over {(px)^n}},{{b}\over {(py)^n}},{{c}\over {(xy)^n}}\right)$$
such that the difference between this element and our original one is in the image of $d_1$.
Looking at the third component, we can thus find an integer $m$ and $e,f\in R$ with
$${r\over {x^ky^k}}-{{c}\over {(xy)^n}}={e\over {x^m}}-{f\over {y^m}}.$$
Multiplying this equation by $x^{k+m}y^{k+m}$, we obtain
$$x^my^m(r-cx^{k-n}y^{k-n})\in (x^{k+m},y^{k+m}).$$
Since the map induced by multiplication by $x^my^m$ from $R/(x^k,y^k)$ to $R/(x^{k+m}y^{k+m})$
is injective, we can thus conclude that the image of $r$ in $R/(x^k,y^k)=Q_k$ is
equal to $x^{k-n}y^{k-n}$ times the image of $c$ in $R/(x^n,y^n)=Q_n$.  Thus the map
induced by multiplication by $x^{k-n}y^{k-n}$ from $Q_n$ to $Q_k$ is surjective.

The proof  of the main theorem of this paper proceeds by repeating the process
of integrating as we have outlined in the construction of $F_2$ from
$F_1$ above. In the process, we obtain elements $z_i$ at stage $i$ and
attempt to add a $U^i$ term in such a way to assure that if $F_i(S,U)$ is
the new polynomial, we have $F_i(z,x)\in y^iR$.  
At each stage there are three cases as above.  If the first 
case holds, there is no problem.  If the second case holds, we adjust
$z_i$ and show that we can still continue as before.  If we
are in the third case, 
we have to multiply by a power of $p$ and it is necessary to increase $L$.
 Using lemma \ref{finitelh}, we can show that eventually 
the third case will not occur.  

The other thing we have to keep track of is the divisibility of binomial coefficients
by powers of $p$.  For
this we introduce the function $\tau$ defined as follows.  Let $\overline\tau(n)$ equal
the sum of digits of $n-1$ in its $p$-adic expansion.  Then let $\tau(n)=\overline\tau(n)/(p-1)$.
We let $\tau(1)=0$.  We note that $\tau(n)$ is a rational number with denominator dividing
$p-1$, and we emphasize the fact 
that $\tau(n)$ is defined by the $p$-adic expansion of $n-1$ rather than 
that of $n$.

As an example. we compute $\tau(p^n)$ and $\tau(p^n+1)$. 
Since the $p$-adic expansion of $p^n-1$ consists of $n$ digits each
equal to $p-1$, we have $\tau(p^n)=n(p-1)/(p-1)=n$. The $p$-adic
expansion of $p^n$ consists of a 1 followed by $n$ zeros, so
$\tau(p^n+1)=1/(p-1)$. 

We next compute the difference $\tau(n)-\tau(n-1)$ for any $n>1$.  
Let $a$ be the largest integer for which $p^a$ divides $n-1$ ($a$
will be zero if $p$ does not divide $n-1$). Then the last $a$ digits
of the $p$-adic expansion of $n-2$ are equal to $p-1$ and the  
$a+1$ digit is less than $p-1$.  Thus when we add
1 to $n-2$, the last $a$ digits are replaced with zeros and the
$a+1$ digit is increased by 1.  Hence we have
\begin{equation}
\tau(n)-\tau(n-1)={1\over {p-1}}-a. \label{tau}
\end{equation}

The connection between $\tau$ and 
 binomial coefficients comes from the following lemma.
In this lemma and the following  discussion we use the notation $p^k\| n$ to mean that $p^k$
is the highest power of $p$ which divides $n$.

\begin{lem}\label{bincoef}
Let $0<j\le i\le p^L$ be integers.  Then  $p^k\|\displaystyle{{p^L-j}\choose 
{i-j}}$, where $$k=\tau(j)+\tau(i-j+1)-\tau(i).$$
\end{lem}

{\bf Proof}
We fix $j$ and $L$ and prove the result by induction on $i$.  If $j=i$, then
$\displaystyle{{{p^L-j}\choose
{i-j}}}=1$ and $\tau(j)+\tau(i-j+1)-\tau(i)=\tau(j)+\tau(1)-\tau(j)=0$; since the highest power
of $p$ which divides 1 is $p^0$, the result is true in this case.

Now assume that  $j<i\le p^L$ and 
that the result is true for $i-1$.  We compute how each
side of the equation changes when we go from $i-1$ to $i$.

For the left hand side, we compute
$${{p^L-j}\choose {i-j}} ={{p^L-j}\choose {(i-1)-j}}\times \left[{{p^L-j-(i-j)+1}\over{i-j}}
\right]=$$
$$ ={{p^L-j}\choose {(i-1)-j}}\times \left[{{p^L-i+1}\over{i-j}}
\right].$$
Hence the highest power of $p$ which divides this binomial coefficient remains the same
unless $p$ divides $p^L-i+1$ or 
$p$ divides $i-j$.  Let $p^a\|(p^L-i+1)$ and $p^b\|(i-j)$.  Then
we obtain the highest power of $p$ dividing $\displaystyle{{p^L-j}\choose {i-j}}$ from that for
$\displaystyle{{p^L-j}\choose {(i-1)-j}}$ by adding $a-b$.

We now compute how $\tau(j)+\tau(i-j+1)-\tau(i)$ changes as we pass from $i-1$ to $i$.
We have
$$[\tau(j)+\tau(i-j+1)-\tau(i)]
-[\tau(j)+\tau((i-1)-j+1)-\tau(i-1)]$$
$$=[\tau(j)+\tau(i-j+1)-\tau(i)]
-[\tau(j)+\tau(i-j)-\tau(i-1)]$$
$$=[\tau(i-j+1)-\tau(i-j)]-[\tau(i)-\tau(i-1)].$$ 

We will now express this difference using the 
$a$ and $b$ defined earlier in the 
proof. We have defined $b$ by the condition $p^b\| (i-j)$ and $a$
by the condition $p^a\| (p^L-i+1)$.  However, since $0<j<i\le p^L$,
we have $0<i-1<p^L$, and thus we also have $p^a\| (i-1)$.
Hence from equation \ref{tau} we obtain
$$[\tau(i-j+1)-\tau(i-j)]-[\tau(i)-\tau(i-1)]
$$
$$= {1\over {p-1}}-b-\left({1\over {p-1}}-a\right) = a-b.$$

This is the same quantity as we computed for
the left hand side of the equation, so this proves the
lemma by induction.

\bigskip
The above proof of Lemma \ref{bincoef} is the shortest way to obtain
the result we need for binomial coefficients that appear in Taylor
expansions but it may appear somewhat unmotivated.  We outline
the steps of a more natural approach to this result.

We let $\sigma_p$ denote the function given by $\sigma_p(n)=$ the sum of
the digits in the $p$-adic expansion of $n$ for an integer $n\ge 0$; note
that we thus have $\tau(n)=\sigma_p(n-1)/(p-1)$.

\begin{enumerate}
\item Show by induction that $p^a\|n!$ where $a=(n-\sigma_p(n))/(p-1)$.
\item Use step 1 to derive a formula for the highest power of $p$
that divides the binomial coefficient $\displaystyle{n\choose k}$
for integers $0\le k\le n$.
\item Apply the formula from step 2 to 
compute the highest power of $p$ that divides
$\displaystyle{{p^L-j}\choose {i-j}}$ and prove Lemma \ref{bincoef}.
\end{enumerate}

\section{Proof of the Main Theorem.}
\begin{thm}(Heitmann)  Let $R$ be an integrally closed local
domain of mixed characteristic $p$ and of dimension 3.  
Let $p,x,y$ be a system of parameters for $R$.  Assume that $R$ contains an element $\sigma$
such that $\sigma^{p-1}=p$.  Let $z$ be an element of $R$ and $N$ a positive integer
such that $p^Nz\in (x,y)$.  Then for any rational number $\epsilon>0$,
 there exists a finite extension
$R'$ of $R$ in which
$$p^{\epsilon}z\in (x,y)R'.$$
\end{thm}

The main idea, as in the analogous situation 
 explained in the previous section, is to prove the existence of a monic polynomial $f(T)$
satisfying conditions which imply that if $w$ is a root of $f(T)$, 
then $p^{\epsilon}w$ will be integral over $R$, and if we let $v$ be
the element such that $p^{\epsilon}z=p^{\epsilon}wx+p^{\epsilon}vy$, 
then $p^{\epsilon}v$ is also integral over $R$.  We next
prove the main lemma which shows that this can be done.

We remark that the reason for requiring that $R$ contain the element
$\sigma$ is that we will need to define powers $p^r$ of $p$ where $r$
is a rational number of the form $n/(p-1)$ for some integer $n$.  We
can then define $p^r=p^{n/(p-1)}=\sigma^n$.

\begin{lem}\label{main} 
Let $R,p,x,y,$ and $z$ be as in the statement of Theorem 1, and let $K$ be
a positive integer.  Then there exists an integer $L$ and a  polynomial
$$f(T)=T^{p^L}+a_1T^{p^L-1}+\cdots +a_{p^L}$$
satisfying the following conditions:
\begin{enumerate}
\item $a_j\in p^{K-\tau(j)}$ for $j=1,\ldots,p^L.$
\item If we let $F(S,U)=S^{p^L}+a_1S^{p^L-1}U+\cdots +a_{p^L}U^{p^L}$ as in the previous
section, we have
$$F^{(k)}(z,x)\in y^{p^L-k}R[p^{-1}]$$ for $k=0,\ldots,p^L$.
\end{enumerate}
\end{lem}

{\bf Proof.}
We fix a positive integer $K$.

As outlined in the previous section, we construct the polynomial $F(S,U)$ inductively
starting with degree 1.  Throughout the construction we will
make use of the quantity $E_i$ defined by the formula
$$E_i=K-\tau(i)+\tau(2).$$
Note the $E_i$ is a rational number with denominator dividing $p-1$ and that it
can be positive or negative.

We now state the situation we have after the $i$th step precisely.
We assume that we have, for each integer $k$ with $1\le k\le i$,

\begin{enumerate}
\item A homogeneous
polynomial $F_k(S,U)$ of degree $k$,
\item an integer $L_k>K$ with $p^{L_k}>k$ for all
$k<i$ (we may have $p^{L_i}=i$; in that case $L_i$ will be the $L$ of the
Lemma and we will be done),
\item an element $z_k$ in $R$ with $p^{N_k}z_k\in (x^k,y^k)$ for
some integer $N_k$.
\end{enumerate}
Furthermore, these polynomials, elements of $R$, and integers
will satisfy the following conditions:

\renewcommand{\theenumi}{\Alph{enumi}}

\begin{enumerate}
\item The coefficient of $S^k$ in $F_k(S,U)$ is $
\displaystyle{{p^{L_k}}\choose k}$. Thus
we can write
$$F_k(S,U)={{p^{L_k}}\choose k}S^k+{{p^{L_k}-1}\choose k-1}a_{k1}S^{k-1}U+\ldots +
{{p^{L_k-k}}\choose 0}a_{kk}U^k,$$
where the $a_{kj}$ are elements of $R[p^{-1}]$ but may not be in $R$. 
\item For $j=1,\ldots,k$ we have $a_{kj}\in p^{K-\tau(j)}$.
\item $F_k^{(m)}(z,x)\in y^{k-m}R[p^{-1}]$ for $m=0,\ldots,k$. 
\item If we let $G_1(S,U)=S$
and $G_k(S,U)=p^{-E_k}\mbox{Int}_{p^{L_{k-1}}} F_{k-1}(S,U)$ for $2\le k\le i$,
 then $G_i(z,x)=z_i$.
\end{enumerate}

This is a lengthy induction hypothesis; however, the
induction step is itself quite complicated and this information
from previous steps is used.

We now do the first step, where $i=1$.
  By hypothesis, there is nonnegative
integer $N_1>0$ such that 
$p^{N_1}z\in(x,y)$, so we can find an $a\in R$ such that $p^{N_1}z+ax\in yR$.  Multiplying this expression
by $p^K$ we have that $p^{N_1+K}z+p^Ka \in yR$.
We let $L_1=N_1+K$ and $a_{11}=p^Ka$, and we define
$$F_1(S,U)=p^{N_1+K}S+p^KaU = p^{L_1}S+a_{11}U.$$
Then $F_1$ is a homogeneous polynomial of degree 1
as required. Let $z_1=z$.  Then we have $L_1>K$, $p^{L_1}>1$, 
 and $p^{N_1}z_1\in (x,y)$,
so conditions 1 through 3 are satisfied.

The coefficient of $S$
in $F_1$ is $p^{L_1}= \displaystyle{{p^{L_1}}\choose 1}$, so condition A
holds. 
Since $\displaystyle{{p^{L_i}\choose 0}}=1$, we have
$$F_1(S,U)={{p^{L_1}}\choose 1}S+{{p^{L_1}}\choose 0}a_{11}U.$$
Condition B states that $a_{11}\in p^{K-\tau(1)}R=p^KR$, which is true since we 
defined $a_{11}=p^Ka$.   Condition C
states that $$F_1(z,x)=p^{L_1}z+a_1x=p^{N+K}z+p^Kax\in yR[p^{-1}],$$ 
which is also true by construction (in fact, we have
$p^{N+K}z+p^Kax\in yR$ in this case).
Finally, since $G_1(S,U)=S$, we have $z_1=z=G_1(z,x)$. 
Hence all the necessary conditions are satisfied for $i=1$.  

We now suppose that $L_k$, $z_k$, and $F_k(S,U)$ have
been defined for $1\le k<i$ and define $F_i(S,U)$,
$L_i$, and $z_i$.

As stated in the previous section, the main idea is to integrate $F_{i-1}(S,U)$ and
modify it so that it will satisfy the conditions listed above. 
However,
to make the construction work it is also necessary to multiply the
integral by a power of $p$ using the number $E_i=K-\tau(i)+\tau(2)$
defined above.  

As a first approximation to $F_i$, we let
$$G_i(S,U)=p^{-E_i}\mbox{Int}_{p^{L_{i-1}}}(F_{i-1}(S,U)).$$
(This is the same polynomial as in condition D for $i\ge 2$.)

Let $$z_i=G_i(z,x).$$

We claim that $z_i$ is an element of $R$. To simplify
notation we let $a_j=a_{i-1,j}$ for $j=0,\ldots,i-1$ (with $a_0=1$), 
so that we have 
$$F_{i-1}(S,U)=\sum_{j=0}^{i-1}{{p^{L_{i-1}}-j}
\choose{i-1-j}}a_jS^{i-1-j}U^j.$$ 
We then have that $a_j\in p^{K-\tau(j)}$ for $j=1,\ldots, i-1$ by induction.

By the definitions of $z_i$ and  $G_i(S,U)$, we have that
$$z_i=p^{-E_i}\sum_{j=0}^{i-1}{{p^{L_{i-1}}-j}\choose{i-j}}a_jz^{i-j}x^j.$$
  We show that each term of this sum is
in $R$; to do so we distinguish two cases.  

If
$j=0$, then $a_j=1$, and we must show that $p^{-E_i}
\displaystyle{{p^{L_{i-1}}}\choose{i}}\in R$.
To see this, we note that
the highest power of $p$ dividing $\displaystyle{{p^{L_{i-1}}}\choose{i}}$
is $p^{{L_{i-1}}-k}$,  where
$p^k\|i$ (of course, $k$ can be zero).  
Thus to show that $p^{-E_i} 
\displaystyle{{p^{L_{i-1}}}\choose{i}}\in R$ 
we must show that $$-E_i+L_i-k=-K+\tau(i)-\tau(2)+L_{i-1}-k\ge 0.$$ 
Since $L_{i-1}> K$, this reduces to showing that $\tau(i)+1\ge \tau(2)+k$.  
Since
$p^k$ divides $i$, the $p$-adic expansion of $i-1$ 
ends with $k$ digits equal to $p-1$. Hence $$\tau(i)+1\ge k+1\ge  k+(1/(p-1))=k+\tau(2),$$
so the inequality holds.

We now consider the general case, where $i\le
p^{L_{i-1}}$ and $0<j<i$.  In this case we
use Lemma \ref{bincoef}.   We must show that  $p^{-E_i}
\displaystyle{{{p^{L_{i-1}}-j}\choose{i-j}}}a_j\in R$. Putting together the powers
of $p$ dividing the factors in this product, 
and using the fact that the highest power of $p$ that divides 
$\displaystyle{{p^{L_{i-1}}-j}\choose{i-j}}$ is $\tau(j)+\tau(i-j+1)-\tau(i)$, the inequality to be 
proven is
$$-K+\tau(i)-\tau(2)+\tau(j)+\tau(i-j+1)-\tau(i)+K-\tau(j)\ge 0.$$
This expression simplifies to
$$\tau(i-j+1)-\tau(2)\ge 0,$$
which is
true since $i>j$, so that $i-j+1\ge 2$ and thus $\tau(i-j+1)\ge\tau(2)$.

Thus we have shown that $z_i\in R$.  
We claim that we also have that $p^{N_i}z_i\in (x^i,y^i)R$ for some $N_i$.
If we can show that $z_i\in (x^i,y^i)R[p^{-1}]$, then we can conclude
that $p^{N_i}z_i\in (x^i,y^i)R$ for some $N_i$ by clearing denominators.
To show that $z_i\in (x^i,y^i)R[p^{-1}]$
we use the fact that $G_i(S,U)$ is,  up to a constant which is a unit in $R[p^{-1}]$,
the integral of $F_{i-1}(S,U)$, and
$z_i=G_i(z,x)$.  Hence, again up to constants which again are units
in $R[p^{-1}]$, we have that 
$G_i^{(k)}(z,x)$ agrees with $F_{i-1}^{(k-1)}(z,x)$
for $k=1,\ldots,i$.  
By induction we have that $F_{i-1}^{(k-1)}(z,x)\in y^{(i-1)-(k-1)})R[p^{-1}]$ for
$k=1,\ldots,i$. It thus follows that $G_i^{(k)}(z,x)\in y^{i-k}R[p^{-1}]$
for $k=1,\ldots,i$.  Hence Lemma \ref{inxy}  implies that $G_i(z,x)\in
(x^i,y^i)R[p^{-1}]$ as was to be shown.  

To summarize the argument up to this point, we have defined 
$$G_i(S,U)=p^{-E_i}\mbox{Int}_{p^{L_{i-1}}}(F_{i-1}(S,U)),$$
defined $z_i=G_i(z,x)$, and shown that $z_i$ is an element of $R$ such
that $p^{N_i}z_i\in(x^i,y^i)$ for some integer $N_i$. 
We now wish to add a $U^i$  term to $G_i(S,U)$ so that
if we evaluate the resulting polynomial at $(z,x)$  
the result is in $y^iR[p^{-1}]$ and
continue the induction. In fact, $G_i(S,U)$ may have to
be modified first as we show below.   We distinguish three cases.

{\bf Case 1.}  The first, and easiest, case, is when $z_i\in (x^i,y^i)R$. 
 We then  know that  
there exists an element $a$ such that 
$$z_i+ax^i\in y^iR.$$
We now let $$F_i(S,U)=p^{E_i}(G_i(S,U)+aU^i).$$
 It
is clear that $F_i(z,x)=p^{E_i}(z_i+ax^i)\in y^iR[p^{-1}]$.  Furthermore, since
the derivatives of $F_i$ with respect to $S$ are, up to a constant which
is a unit in $R[p^{-1}]$, the same as those of $G_i$, it follows from the 
argument in the next
to the last paragraph that
$$F_i^{(k)}(z,x)\in y^{i-k}R[p^{-1}]$$
for $k=1,\ldots,i$, and hence, since this condition also holds for
$k=0$,  it holds for $k=0,\ldots,i$.

We let $L_i=L_{i-1}$.  Then $L_i>K$, and conditions (1) through
(3) are satisfied except that we may have $p^{L_i}=i$, which would
mean that we could not continue the induction to $i+1$.  However,
we claim that in this case the polynomial
$$f(T)=F_i(T,1)$$
satisfies the conditions required in the conclusion of Lemma \ref{main}
and we are done. However, we first verify that conditions A through D
hold.

Conditions A and D hold by construction, and we have already checked condition
C.  That leaves condition B, that the coefficients $a_{ij}$ are in
$p^{K-\tau(j)}R$.  For $j<i$, we have $a_{ij}=a_{i-1,j}$, and
this condition follows by induction.
 Thus we
must check that $a_i\in p^{K-\tau(i)}R$.
Since $a\in R$, we have $a_i=p^{E_i}a\in
p^{E_i}R=p^{K-\tau(i)+\tau(2)}R\subseteq p^{K-\tau(i)}R$, which proves 
that the condition holds in this case as well.

We now return to the case in which $p^{L_i}=i$.  If we
let $f(T)=F_i(T,1)$, then
the leading coefficient of $f(T)$  is $\displaystyle{{p^{L_i}}\choose i}$,
which is 1, so $f(T)$ is monic as required. 
Furthermore, conditions (1) and (2) of Lemma \ref{main} follow immediately
from conditions B and C of the induction hypothesis.
Thus if $p^{L_i}=i$, the proof of Lemma \ref{main} is complete.

If $p^{L_i}>i$, we continue the induction. 
We remark that the argument that $a_{ii}\in p^{K-\tau(i)}R$ and the
fact that the proof of Lemma \ref{main} is complete when $p^{L_i}=i$ 
will be used in the other cases as well.

{\bf Case 2.}  We assume now that $z_i$ is not in $(x^i,y^i)$, but that 
we can write $z_i$ modulo $(x^i,y^i)$ in terms of $z_k$ for $1\le k<i$ in the
following sense: there exist $c_1,\ldots, c_{i-1}$ such that
$$z_i+c_{i-1}xyz_{i-1}+\cdots +c_1x^{i-1}y^{i-1}z_1\in (x^i,y^i)R.$$
We recall that for each $k=1,\ldots i-1$ we have
$$z_k=G_k(z,x),$$
where $G_k$  is the homogeneous polynomial of degree $k$ defined above.
We now let
$$\tilde F_i(S,U)=G_i(S,U)+c_{i-1}yUG_{i-1}(S,U)+\cdots +c_1y^{i-1}U^{i-1}G_1(S,U).$$
From the above equation we thus have
$$\tilde F_i(z,x)\in (x^i,y^i),$$
so we can find an element $a\in R$ such that 
$$\tilde F_i(z,x)+ax^i\in y^iR.$$
We define
$$F_i(S,U)=p^{E_i}(\tilde F_i(S,U)+aU^i).$$

We must now check that $F_i$ has the required properties. 
We let $L_i=L_{i-1}$ as in the previous case.
Since the coefficient of $S^i$ in $F_i$ is the same as that of 
$p^{E_i}G_i=\mbox{Int}_{p^{L_{i-1}}}(F_{i-1})$, we have by induction that this
coefficient is $\displaystyle{p^{L_i}\choose i}$ as required.  
We now check that $a_{ij}\in 
p^{K-\tau(j)}R$ for $j=1,\ldots,i-1$; the verification that $a_{ii}\in
p^{K-\tau(i)}R$ is the same as in the previous case.  We must show that
the coefficients of $p^{E_i}\tilde F(S,U)$ satisfy Condition B,
and, since the coefficients of $p^{E_i}G_i(S,U)$ satisfy this condition,
it suffices to show that the contributions of the 
coefficients of $p^{E_i}y^{i-k}U^{i-k}G_k$
satisfy the condition for each $k=1,\ldots,i-1$.

Let $k$ be an integer with $1\le k \le i-1$.  By construction, we
can write 
$$G_k(S,U)=p^{-E_k}\left(\sum_{j=0}^{k-1}{{p^{L_k}-j}\choose k-j}a_{kj}S^{k-j}U^j\right),$$
where $a_{kj}\in p^{K-\tau(j)}R.$  What we have to show is that if we write
$$p^{E_i}U^{i-k}G_k(S,U)=\sum_{m=0}^{i-1}{{p^{L_i}-m}\choose {i-m}}a'_mS^{i-m}U^m,$$
then we have $a'_m\in p^{K-\tau(m)}R.$

We note that the since we are multiplying by $U^{i-k}$, the term of
$G_k(S,U)$ which contains $a_{kj}$ corresponds to the term of
$p^{E_i}U^{i-k}G_k(S,U)$ which contains $a'_{j+i-k}$.
Thus we have to show that if $a_{kj}\in p^{K-\tau(j)}R$ and
$$p^{E_i}p^{-E_k}{{p^{L_k}-j}\choose k-j}a_{kj}={{p^{L_i}-(j+i-k)}\choose {i-(j+i-k)}}a'_{j+i-k},$$
then $a'_{j+i-k}\in p^{K-\tau(j+k-i)}R$. Dividing by the binomial coefficient 
on the right, this means that 
$$p^{E_i}p^{-E_k}{{p^{L_k}-j}\choose k-j}p^{K-\tau(j)}{{p^{L_i}-(j+i-k)}\choose {i-(j+i-k)}}^{-1}
\in p^{K-\tau(j+i-k)}R.$$

To verify this statement we  use Lemma \ref{bincoef}
to determine the power of $p$  
dividing the binomial coefficients in this expression.  
We also use the above expressions for $E_i$ and $E_k$. What
results is that we must prove the inequality
$$(K-\tau(i)+\tau(2))-(K-\tau(k)+\tau(2))+(\tau(j)+\tau(k-j+1)
-\tau(k))$$
$$+(K-\tau(j))-
(\tau(j+i-k)+\tau(k-j+1)-\tau(i))\ge 
 K-\tau(j+i-k).$$
When this is worked out, it is seen that the two sides of the equation are in fact equal.

We must now prove the Condition C on derivatives.  By construction, we have
$$F_i(z,x)\in y^iR[p^{-1}].$$
To prove that $F_i^{(k)}(z,x)\in y^{i-k}R[p^{-1}]$ for $k\ge 1$ it suffices to
prove the corresponding condition for $\tilde F_i$, and to prove this it suffices to
prove it for each term in the sum defining $\tilde F_i$.  For $G_i$ the proof
is the same as in the previous case.  For $c_{m}y^{i-m}U^{i-m}G_{m}(S,U)$,
it suffices to show that if we let $H(S,U)=y^{i-m}G_m(S,U)$, then $H^{(k)}(z,x)\in
y^{i-k}R[p^{-1}]$ for $k\ge 0$.  By construction, we have that
$G_m^{(k)}(z,x)\in y^{m-k}R[p^{-1}]$ for $k\ge 0$.  Hence
$$H^{(k)}(z,x)=y^{i-m}G_m^{(k)}(z,x)\in y^{i-m}y^{m-k}R[p^{-1}]=y^{i-k}R[p^{-1}]$$
as was to be shown.  

This finishes the second case.

{\bf Case 3.}  This is the case in 
which $z_i$ cannot be forced into $(x^i,y^i)R$ even 
after modifying it by multiples of  $x^{i-k}y^{i-k}z_k$.  
Let $N_i$ be a positive integer
such that $p^{N_i}z_i\in (x^i,y^i)R$.  We now let $L_i=L_{i-1}+N_i$. Note that
this is the only case where we increase $L_i$.   We now let
$$\tilde F_i(S,U)=
\displaystyle{{p^{L_i}\choose i}}{p^{L_{i-1}}\choose i}^{-1}G_i(S,U).$$
Regardless of what $i$ is, the highest power of $p$ dividing 
$\displaystyle{{p^{L_i}\choose i}}{p^{L_{i-1}}\choose i}^{-1}$ 
is $p^{L_i-L_{i-1}}=p^{N_i}$.
Hence, since $p^{N_i}z_i\in (x^i,y^i)R$, we have $\tilde F_i(z,x)\in (x^i,y^i)$,
and we can find an element $a$ such that $\tilde F_i(z,x)+ax^i\in y^iR$.
We now let
$$F_i(S,U)=p^{E_i}(\tilde F_i(S,U)+aU^i).$$

Since $F_i$ is, apart from the 
$U^i$ term, a constant multiple of $G_i$ by a unit
in $R[p^{-1}]$, the condition on derivatives is clear in this case. Furthermore,
we have constructed $F_i$ so that the coefficient of $S^i$ is
$\displaystyle{{p^{L_i}\choose i}}$, so it satisfies the Condition A.
We now check the Condition B on the $a_{ij}$.
In fact, it follows 
from Lemma \ref{bincoef} that the highest power
of $p$ dividing $\displaystyle{{p^{L_i}-j\choose i-j}}$ 
is the same as the highest power of $p$ dividing
$\displaystyle{p^{L_{i-1}}-j\choose i-j}$ for all $j=1,\ldots, i-1$.  
Hence the divisibility conditions on
$a_{ij}$ for $G_i$ are the same as those for $F_i$. Thus,
since the $a_{ij}$ are, up to units in $R$, multiples of
$a_{i-1,j}$ by $p^{N_i}$, the condition follows.
Thus all the conditions hold, and this completes the proof of Case 3.

To finish the proof of the lemma, we invoke Lemma \ref{finitelh}.
As in that lemma, we let $Q_i$ denote the submodule of
$R/(x^i,y^i)$ consisting of elements $u$ annihilated
by $p^N$ for some $N$.  Lemma \ref{finitelh} states that for large $i$
multiplication by $x^ky^k$ identifies $Q_i$ with $Q_{i+k}$
for all $k\ge 0$.  Let $i$ be large enough so that this
holds.
For each $n$ let $M_n$ denote the submodule of $Q_n$ generated by $z_n,
xyz_{n-1},\ldots,x^{n-1}y^{n-1}z_1$, and for each $k\ge 0$
let $N_k$ be the submodule of $Q_i$ corresponding to $M_{i+k}$
under the above identification.  Since the $N_k$ form an
increasing sequence of submodules of $Q_i$, we
must have $N_k=N_{k+1}$ for sufficiently large $k$.  It follows that
multiplication by $xy$ identifies $M_i$ with $M_{i+1}$ for
large $i$.  Thus for large $i$ it is always possible
to write $z_i$ in terms of $z_{i-1},\ldots z_1$ as in
Case 2, and Case 3 will not occur. 
Thus we will eventually reach the situation where
$i=p^{L_i}$, and at that point $F_i(S,U)$ will satisfy the required conditions. 

\bigskip
We now prove the theorem.  Assume that $p^Nz\in (x,y)$, where $p,x,y$ form a
system of parameters.  We wish to show that for all rational $\epsilon>0$, we
have $p^{\epsilon}z\in (x,y)R'$ for some finite extension $R'$ of $R$.  Choose
$K>0$ such that $1/p^K<\epsilon$; we can in fact assume that 
$\epsilon=1/p^K.$   
Let $f(T)$ be a polynomial such that 
\renewcommand{\theenumi}{\arabic{enumi}}
\begin{enumerate} 
\item $f(T)=T^n+a_1T^{n-1}+\ldots +a_n$, where
\item $F^{(k)}(z,x)\in y^{n-k}R[p^{-1}]$ for $k=0,\ldots, n$, and
\item $a_j\in p^{K-\tau(j)}R$ for $j=1,\ldots,n$.
\end{enumerate}

We claim that if we let $g(T)=\sum_{i=0}^np^{{\epsilon}j}a_jT^{n-j}$, then $g(T)$ has
coefficients in $R$.  This amounts to the statement that
$$j\epsilon+K-\tau(j)=j/p^K+K-\tau(j)\ge 0$$ for all $j$.    

To prove this inequality, we first note that
if $j\le p^K$, then $\tau(j)\le K$ and the inequality is clear.
Thus we may assume that $j>p^K$, and this implies that
there is a nonnegative integer $M$ such that $p^{K+M}\le j\le p^{K+M+1}$.

The inequalities $p^{K+M}\le j\le p^{K+M+1}$ imply that
$j/p^K\ge p^M$ and $\tau(j)\le K+M+1$.
We thus have
$$j/p^K+K-\tau(j)\ge p^M+K-(K+M+1)=p^M-M-1.$$
It thus suffices to show that for all nonnegative integers
$M$ and all prime numbers $p$ we have $p^M\ge M+1$.
For $M=0$ this states that $p^0\ge 1$, and for $M=1$ it states
that $p\ge 2$, and both of these statements are true.  The general
case can be shown by an easy induction on $M$.

Thus $g(T)$ has coefficients in $R$.  Let $w$ be a root of
$f(T)$; we then have that $p^{\epsilon}w$ is a root of $g(T)$ so
is integral over $R$.  The condition on derivatives in Lemma \ref{main}
implies that if we let $v=(z-wx)/y$, then $v$ is integral over
$R[p^{-1}]$, and hence $p^{\epsilon}v$ is also integral
over $R[p^{-1}]$.  Furthermore, $p^{\epsilon}v=(p^{\epsilon}z-
(p^{\epsilon}w)x)/y$ is a quotient of an element integral over
$R$ by $y$ so is clearly integral over $R[y^{-1}]$. Thus
$p^{\epsilon}v$ is integral over $R[p^{-1}]$ and $R[y^{-1}]$,
so, since $p$ and $y$ generate an ideal of height 2 in $R$ and
$R$ is a normal domain, $p^{\epsilon}v$ is integral over $R$.
We now have
$$p^{\epsilon}z=(p^{\epsilon}w)x+(p^{\epsilon}v)y,$$
so $p^{\epsilon}z$ is in the ideal generated by $x$ and $y$ in
a finite extension of $R$.
This concludes the proof of the theorem.

\section{A proof of the Canonical Element Conjecture in dimension 3}

In this section we show how Heitmann's theorem can be used to 
give a direct proof of the Canonical
Element Conjecture in dimension 3.  We recall the statement of one version of this
conjecture.

Let $x_1,\ldots, x_d$ be a system of parameters for $R$. Let $K_\bullet$ be the Koszul
complex on $x_1,\ldots, x_d$, and let $F_\bullet $ be a minimal free resolution of
$R/(x_1,\ldots, x_d)$.  Since $K_\bullet$ is a complex of free modules and $F_\bullet$
is exact, there is a map of complexes $\phi_\bullet:K_\bullet\to F_\bullet$ which
is the identity map in degrees $0$ and $1$ (note that the maps $K_1\to K_0$ and
$F_1\to F_0$ are the same).  The Canonical Element Conjecture states that 
$$\phi_d(1)\not\in \m F_d.$$

In dimension 3 in mixed characteristic, we can assume that the system of parameters
is of the form $p^N,x,y$, where $x,y$ form a regular sequence.  Then the homology of
the Koszul complex $K_\bullet$ in degree 1 is isomorphic to $Q/(x,y)$, where 
$Q=\{r\in R|p^Nr\in (x,y)\}$.  If $s\in R$ is such that $sQ\subseteq (x,y)$, then $s$
annihilates the homology of $K_\bullet$ in degree 1, and, since
$K_\bullet$ is exact in degrees 2 and 3, we can construct a map
$\psi_\bullet:F_\bullet
\to K_\bullet$ such that $\psi_\bullet$ is multiplication by $s$
in degrees 0 and 1.  

Now suppose that we had $\phi_d(1)\in \m F_d$.  We take $s=p^e$ for some small rational
number $e$.  Let $R'$ be a finite 
integrally closed extension of $R$ containing $p^e$ and such
that $p^e (QR')\subseteq (x,y)R'$; the existence of such an extension $R'$ is guaranteed by
Heitmann's Theorem.  By the previous paragraph, there is a map 
$\psi_\bullet:F_\bullet\otimes R'\to K_\bullet\otimes R'$
which is multiplication by $p^e$ in degrees 0 and 1. Thus the composition  
$\psi_\bullet \phi_\bullet$ is a map from $K_\bullet\otimes R'$ to itself which is also
given by multiplication by $p^e$ in degrees 0 and 1 (we use the notation
$\phi_\bullet$ also to denote the extension of the original $\phi_\bullet$
to $R'$).  Since $K_\bullet\otimes R'$ is
exact in degrees 2 and 3, we can conclude that multiplication by $p^e$ is homotopic
to $\psi_\bullet \phi_\bullet$. Thus there is a map $\sigma$ from 
$K_2\otimes R'$ to $K_3\otimes R'$
such that for all $a\in K_3\otimes R'$ we have 
$$p^ea = \psi_\bullet \phi_\bullet(a)+\sigma(d^K_3(a)).$$
Since we are assuming that $\phi_d(a)\in \m R'$, and $d^K_3(a)$ is clearly in $\m R'$,
this implies that $p^e\in \m R'$.
The next lemma shows that this is impossible if we take
$e$ small enough.

\begin{lem} Let $R$ be a local integral domain with maximal ideal $\m$,
and let $c$ be a nonzero element of $R$.  Then for
sufficiently small rational $e>0$, for any finite extension $R'$ of $R$, we have
$c^e\not\in\m.$   
\end{lem}
{\bf Proof} We use the fact that there is a discrete valuation $v$ on $R$ such that
$v(m)>0$ for all $m\in \m$. Let $u$ be the minimum of $v(m)$ for a set of generators
$m$ of $\m$.  Then if we choose $e<u/v(c)$, since $v$ can be extended to any finite
extension $R'$ of $R$, we have that $v(c^e)=ev(c)<u\le v(m)$ for all $m$ in a set of
generators for $\m R'$.  Thus $v(c^e)<v(m)$ for all $m\in \m R'$, so $c^e\not\in \m R'$.


\begin{thebibliography}{99}
\baselineskip=20pt

\bibitem{h} {\sc R. Heitmann},
{\it The direct summand conjecture in dimension three}, 
{\it Annals of Mathematics} {\bf 156} (2002), 695--712.
 
\bibitem{b-h} {\sc W. Bruns and J. Herzog}, {\it Cohen-Macaulay Rings},
Cambridge University Press, Cambridge 1993.                                  

\end{thebibliography}
\end{document}